\documentclass{elsarticle}
\usepackage{amssymb, color}
\usepackage{amsmath}
 \usepackage{pifont}
 \usepackage{caption}
 \usepackage[utf8]{inputenc}
\usepackage{setspace}
\makeatletter
\def\LaTeX{\leavevmode L\raise.42ex
    \hbox{\kern-.3em\size{\sf@size}{0pt}\selectfont A}\kern-.15em\TeX}
\makeatother

\newcommand{\BibTeX}{{\rm B\kern-.05em{\sc
          i\kern-.025emb}\kern-.08em\TeX}}

\makeatletter
\def\@currentlabel{2.1}\label{e:dispaa}
\def\@currentlabel{2.21}\label{e:dispau}
\def\@currentlabel{2.22}\label{e:dispav}
\def\@currentlabel{2.23}\label{e:dispaw}
\def\@currentlabel{2.24}\label{e:dispax}
\def\theequation{\thesection.\@arabic\c@equation}
\makeatother

\renewcommand{\theequation}{\arabic{section}.\arabic{equation}}

\newcommand{\R}{\mathbb R}

\def \D{\Delta}

\newtheorem{thm}{Theorem} [section]
\newtheorem{lem}{Lemma} [section]

\newtheorem{cor}{Corollary} [section]
\newtheorem{definition}{Definition} [section]

\newtheorem{rem}{Remark}[section]
\newtheorem{taggedtheoremx}{Theorem}
\newenvironment{taggedtheorem}[1]
 {\renewcommand\thetaggedtheoremx{#1}\taggedtheoremx}
 {\endtaggedtheoremx}

\renewcommand{\theequation}{\thesection.\arabic{equation}}
\renewcommand{\thesection}{\arabic{section}}
\renewcommand{\theequation}{\thesection.\arabic{equation}}
\let\ssection=\section\renewcommand{\section}{\setcounter{equation}{0}\ssection}
\begin{document}

\begin{frontmatter}
\title{Liouville type theorems for stable solutions of  the weighted system involving the
Grushin operator  with negative exponents.}
\author[mf2]{Foued Mtiri}
\ead{mtirifoued@yahoo.fr}
\address[mf2]{ University of Tunis El Manar, Faculty of Sciences of Tunis Department of Mathematics Campus University 2092 Tunis, Tunisia.}

\begin{abstract}
The aim of this paper is to study the stability of solutions to the general weighted system with negative exponents:
\begin{align*}
\Delta_{s} u = \rho(\mathbf{x}) v^{-p}, \quad \Delta_{s} v = \rho(\mathbf{x}) u^{-\theta}, \quad u,v>0 \quad \text{in } \mathbb{R}^N, \quad p \geq \theta > 1, \quad s \geq 0,
\end{align*}
where $\Delta_{s} u = \Delta_{x} u + |x|^{2s} \Delta_{y} u$ is the Grushin operator, and $\rho$ is a nonnegative continuous function satisfying certain conditions.

We show that the system has no stable solution if $p \geq \theta > 1$ and $N_s < 2 \left[ 1 + (2 + \alpha) x_0 \right]$, where $x_0$ is the largest root of the equation:
\begin{align*}
x^4 - \frac{16p\theta(p-1)}{\theta-1} \left( \frac{1}{p+\theta+2} \right)^2 \left[ x^2 + \frac{p+\theta-2}{(p+\theta+2)(\theta-1)} x + \frac{p-1}{(\theta-1)(p+\theta+2)^2} \right]
\end{align*}

This improves previous work in \cite{at}. Our results can also be applied to the weighted equation with negative exponents:
\begin{align*}
\Delta_{s} u = \rho(\mathbf{x}) u^{-p} \quad \text{in } \mathbb{R}^N, \quad \text{where } p > 1.
\end{align*}
\end{abstract}

\begin{keyword}
Stable solutions \sep Liouville-type theorem \sep Grushin operator \sep Critical exponents \sep Elliptic system.
\end{keyword}
\end{frontmatter}

\section{Introduction}
\setcounter{equation}{0}
\medskip

In this paper, we are interested in the classification of stable solutions to the following system:
\begin{align}\label{1.1}
\Delta_{s} u = \rho(\mathbf{x}) v^{-p}, \quad \Delta_{s} v = \rho(\mathbf{x}) u^{-\theta}, \quad u, v > 0 \quad \text{in } \mathbb{R}^N,
\end{align}
where $\Delta_{s} u = \Delta_{x} u + |x|^{2s} \Delta_{y} u$ is the Grushin operator, $s \geq 0$, $1 < \theta \leq p$, and $\rho: \mathbb{R}^N \to \mathbb{R}$ is a radial continuous function satisfying the following assumption:
\begin{enumerate}
  \item[$(\star)$] There exists $\alpha \geq 0$ and $C > 0$ such that $\rho(\mathbf{x}) \geq C \rho_0(\mathbf{x})$ in $\mathbb{R}^N$, where $\rho_0$ is given by
  \[
  \rho_0 := \left( 1 + \|\mathbf{x}\|^{2(s+1)} \right)^{\frac{\alpha}{2(s+1)}},
  \]
  and
  \[
  \|\mathbf{x}\| = \left( |x|^{2(s+1)} + |y|^2 \right)^{\frac{1}{2(s+1)}}, \quad s \geq 0, \quad \mathbf{x} = (x, y) \in \mathbb{R}^N = \mathbb{R}^{N_1} \times \mathbb{R}^{N_2},
  \]
  is the norm corresponding to the Grushin distance, where $|x|$ and $|y|$ are the usual Euclidean norms in $\mathbb{R}^{N_1}$ and $\mathbb{R}^{N_2}$, respectively. The $\|\mathbf{x}\|$-norm is 1-homogeneous for the group of anisotropic dilations associated with $\Delta_{s}$. The operator $\Delta_{s}$ belongs to the class of degenerate elliptic operators, which has received significant attention over the years, see \cite{ DDM, rahal, xxy}. For elementary properties and typical examples of $\Delta_{s}$, we refer to \cite{nv1b, nb}, and we also highlight the recent survey \cite{nvb}.
\end{enumerate}

\medskip

We begin by considering the well-known weighted Lane-Emden system:
\begin{align}\label{1.2kk}
-\Delta u = \rho v^p, \quad -\Delta v = \rho u^\theta, \quad u, v > 0 \quad \text{in } \mathbb{R}^N, \quad p \geq \theta > 1,
\end{align}
where $\rho(x)$ is a radial function satisfying $\rho(x) \geq A (1 + |x|^2)^{\frac{\alpha}{2}}$ at infinity. This system has attracted considerable attention in recent years, see the pioneering works \cite{cow, Hfh, HU, MY}. Cowan \cite{cow} classified positive stable solutions for $N \leq 10$ and $p \geq \theta > 2$, and this result was extended in \cite{HU} for $\rho \equiv (1 + |x|^2)^{\frac{\alpha}{2}}$, with $2 \leq \theta \leq p$ and $N \leq 10 + 4\alpha$. This was further improved by Hajlaoui et al. \cite{Hfh}, who established a new comparison property crucial for dealing with the case $1 < \theta \leq \frac{4}{3}$.

Among other things, Mtiri and Ye \cite{MY} completely classified positive solutions stable outside a compact set for subcritical $(p, \theta)$ pairs.

\medskip

A natural question in the study of systems involving the Grushin operator is whether similar classifications can be obtained as for the Laplace operator.

For the general system with $s \geq 0$:
\begin{align}\label{N1b.1}
-\Delta_{s} u = v^p, \quad -\Delta_{s} v = u^\theta, \quad u, v > 0 \quad \text{in } \mathbb{R}^N = \mathbb{R}^{N_1} \times \mathbb{R}^{N_2}, \quad p \geq \theta > 1,
\end{align}
the Liouville property is less understood and is more complicated to analyze compared to the case when $s = 0$, due to the lack of symmetry in the operator $\Delta_{s}$ and its degeneration on the manifold $\{0\} \times \mathbb{R}^{N_2}$, which introduces additional mathematical difficulties.

Adopting the approach used in \cite{cow, Hfh}, the authors of \cite{DP} extended Cowan's results (Theorem 1 with $s = 0$) and demonstrated that no smooth stable solution exists for \eqref{N1b.1} if $\frac{4}{3} < \theta \leq p$ and $N_s := N_1 + (1+s) N_2$ satisfies:
\[
N_s < 2 + 2\beta_1 t_1^+, \quad \text{where} \quad t_1^\pm = \sqrt{\omega} + \sqrt{\omega - \sqrt{\omega}}, \quad \omega = \frac{p\theta \alpha_1}{\beta_1}, \quad \alpha_1 = \frac{2(p+1)}{p\theta - 1}, \quad \beta_1 = \frac{2(\theta + 1)}{p\theta - 1}.
\]
They also classified bounded stable solutions for $1 < \theta \leq \frac{4}{3}$ and found that $N_s$ satisfies:
\[
N_s < 2 + \left[ 2 + \frac{2(p+1)}{p\theta - 1} + \frac{4(2-p)}{\theta + p - 2} \right] t_1^+.
\]
This result was improved in a work by Mtiri \cite{f}, where the range of nonexistence results was extended beyond the work of \cite{DP}.

\medskip

A new approach, independently developed by Mtiri \cite{ff}, enables a Liouville theorem for stable solutions of \eqref{N1b.1} for any $p, \theta > 0$ and $N_s$ satisfying:
\[
N_s < 2 + \alpha_1 + \beta_1.
\]

\medskip

On the other hand, the nonexistence of stable solutions for systems with negative exponents has attracted significant attention \cite{at}, but the problem remains incompletely addressed.

In the special case where $\rho \equiv 1$, the system \eqref{1.1} becomes:
\begin{align}\label{1b.1}
\Delta_{s} u = v^{-p}, \quad \Delta_{s} v = u^{-\theta}, \quad u, v > 0 \quad \text{in } \mathbb{R}^N, \quad p \geq \theta > 1.
\end{align}
Following the approach in \cite{DP}, Duong-Nguyen \cite{at} established the following nonexistence result:
\begin{taggedtheorem}{A}
  Assume that $1<\theta\leq p$ and
  \[
  N_{s}<2+ \frac{8}{p+\theta+2}t_0^+, \quad t_{0}^{+} = \sqrt{\gamma} + \sqrt{\gamma + \sqrt{\gamma}} \quad \text{and} \quad \gamma = \frac{p\theta(p-1)}{\theta-1},
  \]
  then \eqref{1b.1} has no bounded stable solution.
\end{taggedtheorem}

\smallskip

In this paper, our aim is to generalize \cite{at}. Let us first recall the notion of stability, which is motivated by \cite{Mo}, see also \cite{cow, DP, Hfh}.
\begin{definition}
  A positive solution $(u, v) \in C^2(\mathbb{R}^N) \times C^2(\mathbb{R}^N)$ of \eqref{1.1} is called stable if there are positive smooth functions $\varphi$, $\chi$ such that
  \[
  -\Delta_{s} \varphi = p v^{-p-1}\chi, \quad -\Delta_{s} \chi = \theta u^{-\theta-1}\varphi \quad \text{in} \, \mathbb{R}^N.
  \]
\end{definition}

\medskip

Our main results state as follows:
\begin{thm}
\label{main3}
  Suppose that $\rho$ satisfies $(\star)$ and let $x_0$ be the largest root of the polynomial
  \[
  H(x) = x^{4} - \frac{16p\theta(p-1)}{\theta-1} \left( \frac{1}{p+\theta+2} \right)^2 \left[ x^2 + \frac{p+\theta-2}{(p+\theta+2)(\theta-1)}x + \frac{p-1}{(\theta-1)(p+\theta+2)^{2}} \right].
  \]
  If $1<\theta\leq p$, then \eqref{1.1} has no bounded stable solution if
  \[
  N_{s}<2\left[1 + (2+\alpha) x_0 \right].
  \]
\end{thm}

\medskip

The famous weighted Grushin equation:
\begin{align}\label{Db}
  -\Delta_{s} u = \rho |u|^{p-1} u \quad \text{in} \, \mathbb{R}^N := \mathbb{R}^{N_1} \times \mathbb{R}^{N_2},
\end{align}
has played an important role in the development of nonlinear analysis in the last decades.

\medskip

Firstly, we recall that in the case $s = 0$ and $\rho = 1$, the finite Morse index solutions to the corresponding problem

\begin{align}\label{Db1}
  -\Delta u = |u|^{p-1} u \quad \text{in} \, \mathbb{R}^N, \quad \text{where} \; p > 1,
\end{align}
have been completely classified by Farina \cite{Far}. Indeed, he proved that a smooth nontrivial solution to \eqref{Db1} exists if $1<p< p_{JL}$ and $N \geq 2$. Here, $p_{JL}$ stands for the Joseph-Lundgren exponent (see \cite{Far}) (see also \cite{CW}). Later on, Fazly \cite{Fa} utilized Farina’s approach to obtain the nonexistence for nontrivial stable solutions of
\[
  -\Delta u = |x|^{\alpha} |u|^{p-1} u
\]
when $p>1$ and $N$ satisfying
\[
  N < N_{\alpha} \quad \text{where} \quad N_{\alpha} = 2 + \frac{2(2+\alpha)}{p-1}\left(p + \sqrt{p^2 - p}\right).
\]

\medskip

A large amount of work has been done generalizing this result in various directions. To cite a few, we refer to \cite{DYfa, DY, Hfh, HU, Hfh3}.

\medskip

Let us comment on related results. For semilinear equations with negative exponents
\begin{align}
  \Delta u = \frac{1}{u^p} \quad u > 0 \quad \text{in} \; \mathbb{R}^N, \quad \text{where} \; p > 1,
\end{align}
the finite Morse index solutions have also been classified by Esposito \cite{n2}. See also  \cite{n1, n3, n4, n5}.

\medskip

Another possible generalization corresponds to elliptic problems involving the Grushin operator, i.e. problem \eqref{Db}. In \cite{rahal}, using Farina's approach, Rahal established a Liouville-type theorem for the equation \eqref{Db} with $\rho = |x|_s^{\alpha}$, and $N_s < N_{\alpha}$.

\smallskip

After that, there have been many contributions to the classification of stable solutions to elliptic equations with negative exponents in various cases of nonlinearities \cite{n6, at}.

\medskip

Here, we obtain a classification result for the weighted Grushin equation with negative exponents by studying the system. In fact, when $p = \theta$, using a Souplet-type estimate (see (2.1) below), the system \eqref{1.1} is reduced to the weighted Grushin equation:
\begin{align}\label{D}
 \Delta_{s} u =\rho( \mathbf{x}) u^{-p}, \quad u>0 \quad \mbox{in }\;\; \mathbb{R}^N \quad \mbox{where }\;\;,p>1.
\end{align}

\medskip

As a consequence of Theorem \ref{main3}, we can claim:

\begin{cor}
\label{main2}
  Suppose that $\rho$ satisfies $(\star)$. Let $p > 1$, and
  \[
  N < 2 + \frac{2(2+\alpha)}{p+1} \left( p + \sqrt{p^2 + p} \right).
  \]
  Then \eqref{D} has no bounded stable solution.
\end{cor}

\medskip

\begin{rem}
  \begin{itemize}
    \item If $\alpha = 0$, then the results in Corollary \ref{main2} coincide with those in \cite{Hfh}.
    \item We can show that $\frac{2t_0^+}{p+\theta+2} < x_0$ for any $1 < \theta \leq p$, (see Remark 2.1 below), where $x_0$ is the largest root of the polynomial $H$ given by  Theorem \ref{main3}. Hence, the range of the nonexistence result in Theorem \ref{main3} is larger than that provided by Theorem {\bf A} with $\alpha = 0$.
    \item We note also that the method used in the present paper can be applied to study weighted systems and to a more general class of degenerate operators, such as the $\Delta_s$ operator (see \cite{nb9, nb}) of the form
    \[
    \Delta_s := \sum_{j=1}^N s_j^2 \Delta_{x^{(j)}}, \quad s := (s_1, \dots, s_N): \mathbb{R}^N \to \mathbb{R}^N,
    \]
    where $s_i: \mathbb{R}^N \to \mathbb{R}$ for $i = 1, \dots, N$, are nonnegative continuous functions satisfying some properties such that $\Delta_s$ is homogeneous of degree two with respect to a group dilation in $\mathbb{R}^N$.
  \end{itemize}
\end{rem}

\medskip

This paper is organized as follows. In Section \ref{s3}, we prove comparison properties between $u$ and $v$ of solutions to \eqref{1.1}, and integral estimates derived from the stability. The proof of Theorem \ref{main3} and Corollary \ref{main2} are given in Section \ref{s4}.

\section{Preliminaries}\label{s3}
\setcounter{equation}{0}

In this section, we introduce some notations and prove some important estimates which will be used in this paper.

 \medskip

In the following, $C$ will denote a generic positive constant independent on $(u,v)$, which could be changed from one line to another. The ball of center $0$ and radius $r > 0$ will be denoted by $B_r$.
 \medskip

\subsection{Notation and some known facts about $\D_{s}$}
\medskip
 We now introduce some  notations and basic properties related to the Grushin operator. Let $s$ be a positive real number and let $\mathbf{z}:=(x_{1},\dots,x_{N_{1}},y_{1},\dots,y_{N_{2}})=(x, y)\in \mathbb{R}^{N_1}\times \mathbb{R}^{N_2}=\mathbb{R}^N$  with $N_{1},N_{2}\geq1$ and $N=N_{1}+N_{2}.$  We denote by  $|x|$ (resp. $|y|$ ) the Euclidean norm in $\mathbb{R}^{N_ 1}$ (resp.  $\mathbb{R}^{N_2}$): $|x|:=\sqrt{x_{1}^{2}+\dots+x_{N_{1}}^{2}}$ (resp. $|y|:=\sqrt{y_{1}^{2}+\dots+y_{N_{2}}^{2}}$).
\smallskip

Set $\nabla_x$ and $\nabla_y$ as Euclidean gradients with respect to the variables $x\in \mathbb{R}^{N_1}$  and $y\in \mathbb{R}^{N_2}$ respectively. For $i=1,\dots,N_ 1$ and $j=1,\dots,N_ 2$,  consider the vector fields
$$X_{i}:=\frac{\partial }{\partial x_{i}},\quad\mbox{and} \;\; Y_{i}:=|x|^{s}\frac{\partial}{\partial y_{j}},$$
and the associated gradient as follows
$$\nabla_{s}:=(X_{1}\dots X_{N_{1}},Y_{1}\dots Y_{N_{2}})=(\nabla_x,  |x|^{s}\nabla_y ).$$
The Grushin operator $\Delta_{s}$, is the operator defined by

$$\Delta_{s}:=\nabla_{s}\cdot \nabla_{s}=\sum_{i=1}^{N_{1}}X_{i}^{2}+\sum_{j=1}^{N_{2}}Y_{j}^{2}=\Delta_{x}+|x|^{2s}\Delta_{y},$$
where $\Delta_{x}$ and $\Delta_{y}$ are Laplace operators in the variables $x\in \mathbb{R}^{N_1}$ and $y\in \mathbb{R}^{N_2}$ respectively. The anisotropic dilation attached to $\Delta_{s}$ is given by
$$\delta_\lambda(\mathbf{z})=(\lambda x, \lambda^{1+s}y), \quad \lambda>0 \quad\mbox{and}\quad \mathbf{z}:=(x, y)\in \mathbb{R}^N=\mathbb{R}^{N_1}\times \mathbb{R}^{N_2}.$$
It is easy to see that $\Delta_{s}$ is the homogeneous of degree two with respect to the
dilation, i.e,

$$\Delta_s(\delta_{\lambda}u)=\lambda^{2}\delta_{\lambda}(\Delta_s u),$$
where $\delta_\lambda u(\mathbf{z})=u(\lambda x, \lambda^{1+s}y).$  The change of variable formula for the Lebesgue
measure gives that

$$d \delta_\lambda(\mathbf{z})=\lambda^{N_s}dxdy,$$
where
$$N_s:=N_1+(1+s)N_2,$$
is the homogeneous dimension with respect to dilation $\delta_\lambda$ and $dxdy$ denotes the Lebesgue measure on  $\mathbb{R}^N.$

\medskip

In the following, we show some properties of the polynomials $L$ and $H$, useful for our proofs. Let

 \begin{align}\label{L}
 L(z):=z^4-\frac{16p\theta(p-1)}{\theta-1}z^2-\frac{16p\theta(p-1)(p+\theta-2)}{(\theta-1)^2}z-\frac{16p\theta(p-1)^2}{(\theta-1)^2}.
\end{align}

\medskip
\subsection{Property of the polynomial $H$}

\medskip

\begin{lem}
\label{L} $L(2t_0)<0$ and $L$ has a unique root $z_0$ in the interval
$(2t_0, \infty)$.
\end{lem}

\noindent
{\bf Proof.} A simple computation leads to
$$\frac{t_0^2}{2t_0+1}=\sqrt{\gamma }\quad \mbox{with} \quad t_{0}^{+} = \sqrt{\gamma}+\sqrt{\gamma+\sqrt{\gamma}} \quad \quad\mbox{and} \quad \gamma =\frac{p\theta(p-1)}{\theta-1}$$
 Obviously
\begin{align*}
L(2t_0)=& \;16t_0^4 - 16\gamma\left[4t_0^2+\frac{2(p+\theta-2)}{\theta-1}t_0+\frac{p-1}{\theta-1}\right]\\
=& \;16\gamma(2t_{0}+1)^{2} - 16\gamma\left[4t_0^2+\frac{2(p+\theta-2)}{\theta-1}t_0+\frac{p-1}{\theta-1}\right]\\
=& \;  16\gamma\left[(2t_{0}+1)^{2}-4t_0^2-\frac{2(p+\theta-2)}{\theta-1}t_0-\frac{p-1}{\theta-1}\right]\\
=& \;  16\frac{p\theta(p-1)}{(\theta-1)^{2}}(2t_{0}+1)(\theta-p)
\end{align*}

As $p\geq\theta>1,$ we have then  $L(2t_0)<0$. Now we consider $L(p)$. Rewrite
\begin{align*}
L(z) = z^4 - 16\frac{p\theta(p-1)}{\theta -1}\left(z^2 + \frac{p+\theta-2}{\theta-1}z + \frac{p-1}{\theta-1}\right).
\end{align*}
For $s > 0$, we see that
\begin{align*}
\left(z^2 + \frac{p+\theta-2}{\theta-1}z + \frac{p-1}{\theta-1}\right)_p' = \frac{1}{\theta-1}(z +1) > 0,
\end{align*}
Then for $z > 0$, as $p\geq\theta>1,$  there holds
\begin{align*}
 z^2 + \frac{p+\theta-2}{\theta-1}z + \frac{p-1}{\theta-1}> z^2 + 2z + 1 = (z+1)^2 \quad\mbox{and} \quad \frac{p\theta(p+1)}{\theta + 1} \geq \theta^2.
\end{align*}
Finally, we get (for $\theta> 1$)
\begin{align*}
L(\theta) < \theta^4 - 16\theta^2(\theta-1)^2=-\left(15\theta^{4}+16(2\theta+1)\right)<0,
\end{align*}
and
\begin{align*}
L'(\theta) = 4\theta^3 - 16\frac{p\theta(p-1)}{\theta - 1}\left(2\theta + \frac{p+\theta-2}{\theta-1}\right) < 4\theta^3 - 16\theta^2(2\theta + 2) = -4\theta^{2}\left(7\theta+8)\right)<0,
\end{align*}
Moreover,  we have
$$L''(z)=12z^2-\frac{32p\theta(p-1)}{\theta-1},$$ then
$L''$ can change at most once the sign from negative to positive for $z\geq 2$. As
$\lim_{s\rightarrow\infty}L(z)= \infty$, it's clear that $L$ admits a unique root in $(2t_0^+,\infty).$  Hence, there holds $2t_0^+ < z_0.$ \qed
\medskip

\begin{rem}
Performing the change of variables $x=\frac{1}{p+\theta+2}z$ in \eqref{L}, a direct calculation yields

\begin{align*}
    H(x)&=:\left(\frac{1}{p+\theta+2}\right)^4L(z)\\
    & = x^{^{4}}- \frac{16p\theta(p-1)}{\theta-1}\left(\frac{1}{p+\theta+2}\right)^4\left[z^2+\frac{p+\theta-2}{\theta-1}z+\frac{p-1}{\theta-1}\right]\\
    & =x^{^{4}}- \frac{16p\theta(p-1)}{\theta-1}\left(\frac{1}{p+\theta+2}\right)^2\left[x^2+\frac{p+\theta-2}{(p+\theta+2)(\theta-1)}x+\frac{p-1}{(\theta-1)(p+\theta+2)^{2}}\right]
\end{align*}
 Hence  $ H(x)<0$ if and only if $L(z)<0$.
 \smallskip

Using the above Lemma, $x_0 = \frac{1}{p+\theta+2}z_0$ is the
largest root of the polynomial $H$, and $x_0$ is the unique root of $H$ for $x\geq \frac{2}{p+\theta+2}t_0$.
\end{rem}

\medskip

\subsection{Main technical tools}
\medskip

In order to prove our results, we need some preliminary results for solutions to the system \eqref{1.1},  as integral estimates, comparison property of $u$, $v$ and an integral inequality derived from the stability.

\medskip

The following is a comparison result between the components $u$, $v$ of solutions to the system \eqref{1.1}.

\begin{lem}\label{Soup}(Comparison property.)
 Let $p\geq \theta>1$ and suppose that  $\rho$ satisfies $(\star)$. Let $(u,v)$ be a solution of \eqref{1.1} and assume that $v$ is bounded, then
\begin{equation}\label{estS}
	v^{p-1}\leq \frac{p-1}{\theta-1}u^{\theta-1},
	\end{equation}
and 
  \begin{align}\label{investS}
u\leq \|v\|_\infty^\frac{p-\theta}{\theta-1}v.
 \end{align}

\end{lem}
\noindent
{\bf Proof.} The proof adapt an idea of \cite{DP}, originally coming from \cite{f}. Let  $\sigma=\frac{\theta-1}{p-1}\in (0,1],$ $\lambda = \sigma^\frac{-1}{1-p}$ and $w= v- \lambda u^{\sigma}.$ Since the simple calculation implies that
\begin{align*}
\D_{s} w = \D_{s} v-\lambda\sigma u^{\sigma-1} \D_{s}u -\lambda\sigma (\sigma-1)|\nabla_{s}u|^{2}u^{\sigma-2}&\geq \D_{s} v-\lambda\sigma u^{\sigma-1} \D_{s}u\\
&=\rho( \mathbf{x}) u^{\sigma-1}\left[ u^{-\theta+1-\sigma}- \lambda\sigma v^{-p}\right]\\
 & = \rho( \mathbf{x})u^{\sigma-1}\left[ u^{-\theta+1-\sigma}- \lambda^{p} v^{-p}\right]\\
 & =  \rho( \mathbf{x})u^{\sigma-1}\left[ -\lambda^{p}v^{-p}+u^{-p\sigma} \right] \\
 & = \rho( \mathbf{x})u^{\sigma-1}\left[- \frac {\lambda^{p}}{ v^{p}}+\frac{1}{u^{p\sigma}} \right].
\end{align*}

Therefore, for any $\sigma\in (0,1]$, there exists $C > 0$ such that
\begin{align}\label{xy}
C\frac{u^{\sigma-1}}{v^{p}u^{p\sigma}}\left[v^{p} -\left(\lambda u^{\sigma}\right)^{p} \right] \leq \rho( \mathbf{x})\frac{u^{\sigma-1}}{v^{p}u^{p\sigma}}\left[v^{p} -\left(\lambda u^{\sigma}\right)^{p}\right]\leq \D_{s} w.
\end{align}

We need to prove that $$v\leq\lambda u^{\sigma}.$$
 We shall show that
 \begin{align}\label{2.9}
w \leq0,
\end{align}
 by a contradiction argument. Suppose that
 \begin{align}\label{00xy}
\sup_{\mathbb{R}^N}w>0.
\end{align}
\medskip

Next, we split the proof into two cases.
\medskip

\textit{Case 1: } We consider the case where the supremum of $w$ is attained at infinity.

\medskip
Choose now   $\phi_{R}(x,y)=\psi^{m}(\frac{x}{R},\frac{y}{R^{1+s}}),$  where $m>0,$ and $\psi $ is a cut-off function in $ C_c^\infty\left(\mathbb{R}^N, [0,1]\right),$ such that $$\psi=1 \quad \mbox{on} \quad B_{1}\times B_{1},\quad \mbox{and} \quad  \psi=0 \quad \mbox{outside } \quad B_{2}\times B_{2^{1+s}}.$$

A simple calculation implies that

$$\frac{|\nabla_{s}\phi_{R}|^{2}}{\phi_{R}}\leq  \frac{C}{R^{2}}\phi_{R}^{\frac{m-2}{m}}\quad \mbox{and} \quad  |\D_{s}(\phi_{R})|\leq \frac{C}{R^{2}}\phi_{R}^{\frac{m-2}{m}}.$$

Set $$w_{R}=\phi_{R}w,$$ which is a compactly supported function. Then there exists $(x_{R},y_{R})\in  B_{2R}\times B_{(2R)^{1+s}},$ such that
\begin{align*}
 w_{R}(x_{R},y_{R})=\max_{\mathbb{R}^N}w_{R}(x,y)\rightarrow \sup_{\mathbb{R}^N}w(x,y) \quad \mbox{as }\; R\rightarrow \infty.
\end{align*}
 This implies
\begin{align*}
\nabla_{s} w_{R}(x_{R},y_{R})=0 \quad \mbox{and}\quad \D_{s} w_{R}(x_{R},y_{R})\leq0,
\end{align*}
which means that at $(x_{R},y_{R})$,
\begin{align}\label{2xy}
\nabla_{s} w=-\phi_{R}^{-1}\nabla_{s}\phi_{R}w \quad \mbox{and}\quad \phi_{R}\D_{s} w\leq 2w\phi_{R}^{-1}|\nabla_{s}\phi_{R}|^{2}-w\D_{s}\phi_{R}.
\end{align}
 From \eqref{2xy},  and using the properties of $\phi_{R}$,  we can conclude then

 \begin{align}\label{3xy}
\phi_{R}\D_{s} w\leq \frac{C}{R^{2}}\phi_{R}^{\frac{m-2}{2}}w.
\end{align}

 Furthermore, for  $w= v- \lambda u^{\sigma}\geq0,$ we  observe that
 \begin{align}\label{4xy}
 \frac{v^{p}}{w^{p}}-\frac{(\lambda u^{\sigma})^{p}}{w^{p}}\geq1, \quad \mbox{ or equivalently }\quad
\lambda^{-p}v^{p} - u^{p\sigma}\geq\lambda^{-p} w^{p}.
\end{align}
 Multiplying  \eqref{xy} by $\phi_{R},$ combining  it with \eqref{4xy} and \eqref{3xy}, one obtains

  \begin{align*}
u^{\sigma-1} w^{p}\phi_{R}^{\frac{m+2}{2}}\leq \frac{C}{R^{2}}w\phi_{R}.
\end{align*}
As $\sigma\leq1,$ the sequence $u(x_{R},y_{R})$ is bounded. We choose
 \begin{align*}
\theta=\frac{m+2}{m} \quad \mbox{so that }\;m=\frac{2}{\theta-1},
\end{align*}
there holds then
   \begin{align*}
 w_{R}^{p-1}(x_{R},y_{R})\leq \frac{C}{R^{2}}.
\end{align*}
 Taking the limit $R \to\infty$, we have $\sup_{\mathbb{R}^N}w=0,$ which contradicts \eqref{00xy}, the claim follows.

\medskip
\textit{Case 2: } If there exists $(x^{0}, y^{0}),$ such that $\sup_{\mathbb{R}^N}w=v(x^{0}, y^{0})- \lambda u^{\sigma}(x^{0}, y^{0})>0,$ then

$\frac{\partial w}{\partial x_{i}}(x^{0}, y^{0})=0$ and  $\frac{\partial w}{\partial y_{i}}(x^{0}, y^{0})=0,$ \;\;  $\frac{\partial^{2} w}{\partial x_{i}^{2}}(x^{0}, y^{0})\leq0$ and  $\frac{\partial^{2} w}{\partial y_{i}^{2}}(x^{0}, y^{0})\leq0.$ However, the left-hand side of \eqref{xy} at $(x^{0}, y^{0})$ is positive. Thus we obtain a contradiction. So we are done.

\medskip

To prove \eqref{investS}, consider $w= u- \lambda v,\;$ where $\lambda = \|u\|_{\infty}^\frac{p-\theta}{\theta-1}$ and we will establish again \eqref{2.9}.
 As $p \geq \theta$ and $u$ is bounded, there holds

\begin{align}\label{x1y}
 \begin{split}
\D_{s} w = \rho( \mathbf{x})\left( v^{-p}-\lambda u^{-\theta}\right) &=\rho( \mathbf{x})\left[ -\lambda u^{-\theta}+ \left(\frac{v}{ \|v\|_{\infty}}\right)^{-p} \|v\|_{\infty}^{-p}\right]\\
&\geq\rho( \mathbf{x})\|v\|_{\infty}^{\theta - p} v^{-\theta}-\lambda u^{-\theta}\\
&\geq \rho( \mathbf{x})\|v\|_{\infty}^{\theta - p} \left(-\lambda^{\theta}u^{-\theta}+v^{-\theta}\right).
\end{split}
\end{align}

For the rest of the proof we just replace  \eqref{xy} by  \eqref{x1y}, so we omit the details and the proof is completed.  \qed
\medskip
\subsection{Stability inequality}

\medskip

Inspired by \cite{DP}, we establish the following \textit{a priori} integral estimates for solutions of the Lane-Emden system \eqref{1.1}.
\begin{lem}\label{l.2.2}
If $(u,v)$ is a nonnegative stable solution of \eqref{1.1}.  Then for all $\phi \in C_c^1(\R^N),$ we have
\begin{equation}
\label{1.3}  \sqrt{p\theta}\int_{\R^N} \rho( \mathbf{x})
u^{-\frac{\theta+1}{2}}v^{-\frac{p+1}{2}}\phi^2dxdy  \leq
\int_{\R^N}|\nabla_{s}\phi|^2dxdy ,
\end{equation}
where
$\nabla_s :=(\nabla_x,  |x|^{s}\nabla_y )$ denotes the Grushin gradient.
\end{lem}

\noindent{\bf Proof.} Let $(u,v)$ denote a  stable solution of \eqref{1.1}. By the definition of stability, there exist positive smooth functions $\varphi$, $\psi$ verifying
\begin{align*}
 -\frac{\Delta_{s} \varphi}{ \varphi} = p \rho( \mathbf{x}) v^{-p-1}\frac{\psi}{\varphi}, \quad -\frac{\Delta_{s} \psi}{ \psi} = \theta \rho( \mathbf{x})u^{-\theta-1}\frac{\varphi}{\psi}\quad \mbox{in }\, \R^N.
\end{align*}

Let $\gamma, \chi \in C_c^1(\R^N)$. Multiplying the first equation by $\gamma^{2}$ and the second by $\chi^{2}$ and integrate
over $\R^N,$ we arrive at

\begin{align*}
 p \int_{\R^N}\rho( \mathbf{x})v^{-p-1}\frac{\psi}{\varphi}\gamma^{2}dxdy=
-\int_{\R^N}\frac{\Delta_{s} \varphi}{ \varphi}\gamma^{2}dxdy,
\end{align*}

 and

\begin{align*}
 \theta \int_{\R^N} \rho( \mathbf{x}) u^{-\theta-1}\frac{\varphi}{\psi}\chi^{2}dxdy=
-\int_{\R^N}\frac{\Delta_{s} \psi}{ \psi}\chi^{2}dxdy.
\end{align*}

The simple calculation implies that
\begin{align*}
  \int_{\R^N}\left(-\frac{\Delta_{s} \varphi}{ \varphi}\gamma^{2}-|\nabla_{s}\gamma|^{2}\right)dxdy& = \int_{\mathbb{R}^N}\Big(\nabla_{s}\varphi\cdot \nabla_{s}(\gamma^{2}\varphi^{-1})-|\nabla_{s}\gamma|^{2}\Big)dxdy\\
& = \int_{\mathbb{R}^N}\Big(-\varphi^{-2}|\nabla_{s}\varphi|^{2}\gamma^{2}+2\varphi^{-1}\gamma\nabla_{s}\varphi\cdot\nabla_{s}\gamma-|\nabla_{s}\gamma|^{2}\Big)dxdy
\\
& = \int_{\mathbb{R}^N}-\Big(\varphi^{-1}\gamma\nabla_{s}\varphi-\nabla_{s}\gamma\Big)^{2}dxdy \leq 0
\end{align*}
 Proceeding as above,  we can easily show  that
\begin{align*}
-\int_{\R^N}\frac{\Delta_{s} \psi}{ \psi}\chi^{2}dxdy\leq\int_{\R^N} |\nabla_{s}\chi|^{2}dxdy.
\end{align*}

Using the inequality $2ab\leq a^{2}+b^{2},$   we deduce that

\begin{align*}
2\rho( \mathbf{x}) \sqrt{p\theta v^{-p-1}u^{-\theta-1}\gamma^{2}\chi^{2}}\leq \rho( \mathbf{x}) \Big( pv^{-p-1}\frac{\psi}{\varphi}\gamma^{2}+\theta u^{-\theta-1}\frac{\varphi}{\psi}\chi^{2}\Big).
\end{align*}

Taking $\phi=\chi=\gamma$  and combining all these inequalities,  we get  readily the estimate \eqref{1.3}. \qed

\medskip

At last, using the  stability inequality and  the  comparison property between $u$ and $v$, we can derive the following integral
estimates for all  solutions of the system  \eqref{1.1},  which is crucial for our analysis.

\begin{lem}\label{l.2.1} Let $p\geq \theta> 1$. Suppose that $\rho$ satisfies $(\star).$ then Then, there exists a positive constant
$C>0$ such that for any solution  $(u,v)$  of \eqref{1.1} and $R\geq1,$ there holds

\begin{align}\label{VV}
  \int_{ B_{R}\times B_{R^{1+s}}} v^{-\frac{p+\theta+2}{2}} dxdy\leq C R^{N_{s}-\alpha-2}.
\end{align}
\begin{align}\label{UU}
   \int_{ B_{R}\times B_{R^{1+s}}} u^{-\theta} dxdy \leq CR^{N_{s}-\frac{2(p-1)\theta}{p\theta-1}-\frac{(p\theta-1)\alpha}{\theta(p-1)}}.
\end{align}
\end{lem}
\noindent{\bf Proof.} The proof adapt an idea of \cite{at}. By Lemma 2.2, we conclude that  $$ Cu^{-\frac{(\theta-1)(p+1)}{2(p-1)}}\leq v^{-\frac{p+1}{2}}.$$

Substituting this in \eqref{1.3},
\begin{equation}
\label{1BB.3}  C\int_{\mathbb R^N}\rho( \mathbf{x})u^{\frac{-p\theta+1}{p-1}}\phi^2\leq\sqrt{p\theta}\int_{\R^N} \rho( \mathbf{x})
u^{-\frac{\theta+1}{2}}v^{-\frac{p+1}{2}}\phi^2dxdy  \leq
\int_{\R^N}|\nabla_{s}\phi|^2dxdy .
\end{equation}

Let $\frac{p\theta-1}{p-1}>\theta.$  Applying H\"older's inequality and adding the result to \eqref{1BB.3},  we derive

\begin{align*}
 \int_{\mathbb R^N} u^{-\theta}\phi^2 dxdy
& \;\leq\;C\left( \int_{\mathbb R^N}\rho( \mathbf{x}) u^{\frac{-p\theta+1}{p-1}}\phi^2dxdy\right)^{\frac{\theta (p-1)}{\theta p-1}}
\; \;\;\;\times \left( \int_{\mathbb R^N}\rho( \mathbf{x})^{-\frac{(\theta p-1)^{2}}{\theta( p-1)(\theta -1)}}\phi^2dxdy\right)^{\frac{\theta -1}{\theta p-1}}\\
  & \; \leq  \;C\left( \int_{\mathbb R^N}|\nabla_{s}\phi|^2dxdy\right)^{\frac{\theta (p-1)}{\theta p-1}}
\; \;\;\;\times \left( \int_{\mathbb R^N}\rho( \mathbf{x})^{-\frac{(\theta p-1)^{2}}{\theta( p-1)(\theta -1)}}\phi^2dxdy\right)^{\frac{\theta -1}{\theta p-1}}.
\end{align*}

Let $\chi_{j}\in C_c^\infty\left(\mathbb{R}, [0,1]\right),$ $j=1,2$ be a cut-off function verifying $0 \leq \chi_{j} \leq 1,$ $$\chi_{j}=1 \quad \mbox{on} [-1,1],\quad \mbox{and} \quad  \chi_{j}=0 \quad \mbox{outside } \quad [-2^{1+(j-1)s},2^{1+(j-1)s}].$$
For $R\geq 1,$ put $\psi_{R}(x,y)=\chi_{1}(\frac{x}{R})\chi_{2}(\frac{y}{R^{1+s}}),$ it is easy to verify  that there exists $C >0$ independent of $R$ such that

$$|\nabla_{x}\psi_{R}|\leq \frac{C}{R}\quad \mbox{and} \quad  |\nabla_{y}\psi_{R}|\leq \frac{C}{R^{1+s}},$$

\smallskip

$$|\D_{x}\psi_{R}|\leq \frac{C}{R^{2}}\quad \mbox{and} \quad  |\D_{y}\psi_{R}|\leq \frac{C}{R^{2(1+s)}}.$$

Take $\phi=\psi_{R}^{m},$ hence

\begin{align*}
& \; \int_{\mathbb R^N} u^{-\theta} \psi_{R}^{2m}dxdy\\
& \;
\leq\;CR^{-2\frac{\theta (p-1)}{\theta p-1}}\left( \int_{ B_{2R}\times B_{(2R)^{1+s}}} \psi_{R}^{\frac{2(m-1)}{m}}dxdy\right)^{\frac{\theta (p-1)}{\theta p-1}}
\; \;\;\;\times \left( \int_{ B_{2R}\times B_{(2R)^{1+s}}}\rho( \mathbf{x})^{-\frac{(\theta p-1)^{2}}{\theta( p-1)(\theta -1)}}\psi_{R}^{2m}dxdy\right)^{\frac{\theta -1}{\theta p-1}}\\
  & \; \leq  \; CR^{-2\frac{\theta (p-1)}{\theta p-1}}R^{N_{s}\frac{\theta (p-1)}{\theta p-1}}R^{-\alpha\frac{\theta p-1}{\theta (p-1)}}R^{N_{s}\frac{\theta -1}{\theta p-1}}=CR^{N_{s}-\frac{2(p-1)\theta}{p\theta-1}-\frac{(p\theta-1)\alpha}{\theta(p-1)}},
\end{align*}
which yields the inequality \eqref{UU}. Similarly, we obtain the estimate for $v$. \qed
 \medskip

\section{Proofs of Theorem  \ref{main3} and Corollary \ref{main2}.}\label{s4}
\setcounter{equation}{0}
The following lemma plays an important role in dealing with Theorems \ref{main3} and Corollary \ref{main2}, where we use some ideas from \cite{Hfh, f}. Here and in the following, we define $R_{k} = 2^{k}R$ for all $R > 0$ and integers $k\geq 1$.
\begin{lem}
\label{newl}

\label{newl}Suppose that $\rho$ satisfies $(\star)$ and let  $(u,v)$ be a stable solution of \eqref{1.1}. Then for any $z> \frac{p-1}{2}$ verifying
$L(z) < 0$, there exists $C <\infty$ such that
\begin{equation}
\label{3.1}
\int_{ B_{R}\times B_{R^{1+s}}}  \rho( \mathbf{x})u^{-\theta}v^{-z-1} dxdy \leq\frac{C}{R^2}\int_{ B_{2R}\times B_{(2R)^{1+s}}}v^{-z} dxdy, \quad \forall\; R > 0.
\end{equation}
where $L(z)$  are defined by \eqref{L} .
\end{lem}
\noindent {\bf Proof}. Let $(u,v)$ be a stable solution of
\eqref{1.1}. Let $\phi \in C_0^2(\R^N)$ and $\varphi =
u^{\frac{1-q}{2}}\phi$ with $q > 0$. Integrating by parts, we get
\begin{align}\label{3.2}
\begin{split}
\int_{\mathbb{R}^N}|\nabla_{s}\varphi|^2dxdy
&=\frac{(1-q)^2}{4q}\int_{\mathbb{R}^N} \rho( \mathbf{x})u^{-q}v^{-p}\phi^2dxdy+\int_{\mathbb{R}^N}u^{1-q}|\nabla_{s}\phi|^2dxdy\\
&+\frac{1-q}{4q}\int_{\mathbb{R}^N}u^{1-q}\Delta_{s}(\phi^2)dxdy.
\end{split}
\end{align}
 Take $\varphi$ into the stability inequality \eqref{1.3} and using \eqref{3.2}, we obtain
\begin{align*}
\sqrt{p\theta}\int_{\mathbb{R}^N} \rho( \mathbf{x})u^{-\frac{\theta+1}{2}}v^{-\frac{p+1}{2}}u^{-q+1}\phi^2dxdy&\leq \int_{\mathbb{R}^N}|\nabla_{s}\varphi|^2dxdy\\
&\leq\frac{(-q+1)^2}{4q}\int_{\mathbb{R}^N} \rho( \mathbf{x})u^{-q}v^{-p}\phi^2dxdy\\
&+C\int_{\mathbb{R}^N}u^{-q+1}\Big[|\nabla_{s}\phi|^2+\Delta_{s}(\phi^2)\Big]dxdy,
\end{align*}
so we get
\begin{align*}
a_1\int_{\mathbb{R}^N} \rho( \mathbf{x})u^{-\frac{\theta+1}{2}}v^{-\frac{p+1}{2}}u^{-q+1}\phi^2dxdy\leq
\int_{\mathbb{R}^N} \rho( \mathbf{x})u^{-q}v^{-p}\phi^2dxdy +C\int_{\mathbb{R}^N}u^{-q+1}\Big[|\nabla\phi|^2+\Delta(\phi^2)\Big]dxdy,
\end{align*}
where $a_1=\frac{4q\sqrt{p\theta}}{(-q+1)^2}$. Choose now  $\phi(x,y)=\psi(\frac{x}{R},\frac{y}{R^{1+s}}),$  where $\psi $ a cut-off function in $ C_c^\infty\left(\mathbb{R}^N=\mathbb{R}^{N_1}\times\mathbb{R}^{N_2}, [o,1]\right),$ such that $$\psi=1 \quad \mbox{on} \quad B_{1}\times B_{1},\quad \mbox{and} \quad  \psi=0 \quad \mbox{outside } \quad B_{2}\times B_{2^{1+s}}.$$

A simple calculation,  implies that

$$|\nabla_{s}\phi|\leq \frac{C}{R}\quad \mbox{and} \quad  |\D_{s}(\phi^2)|\leq \frac{C}{R^{2}}.$$

Hence,
 \begin{equation}\label{3.3}
I_1:=\int_{\mathbb{R}^N} \rho( \mathbf{x})u^{-\frac{\theta+1}{2}}v^{-\frac{p+1}{2}}u^{-q+1}\phi^2dxdy\leq
\frac{1}{a_1}
\int_{\mathbb{R}^N} \rho( \mathbf{x})u^{-q}v^{-p}\phi^2dxdy+\frac{C}{R^2}\int_{ B_{2R}\times B_{(2R)^{1+s}}}u^{-q+1}dxdy
\end{equation}

Furthermore,, using $v^{\frac{-r+1}{2}}\phi,$  $r
> 0$ as test function in \eqref{1.3}. As above,  we get readily
\begin{equation}\label{3.4}
I_2 :=\int_{\mathbb{R}^N} \rho( \mathbf{x})u^{-\frac{\theta+1}{2}}v^{-\frac{p+1}{2}}v^{-r+1}\phi^2dxdy\leq
\frac{1}{a_2} \int_{\mathbb{R}^N} \rho( \mathbf{x})u^{-\theta}v^{-r}\phi^2dxdy+\frac{C}{R^2}\int_{ B_{2R}\times B_{(2R)^{1+s}}}v^{-r+1}dxdy.
\end{equation}
 with $a_2=\frac{4r\sqrt{p\theta}}{(-r+1)^2}$. Combining \eqref{3.3} and  \eqref{3.4},  we have then
 \begin{align}
 \label{3.5}
  \begin{split}
 & I_1+{a_2}^\frac{2(r-1)}{p-1} I_2\\
   &\leq \frac{1}{a_1}
\int_{\mathbb{R}^N} \rho( \mathbf{x})u^{-q}v^{-p}\phi^2dxdy+{a_2}^\frac{2r-1-p}{p-1}\int_{\mathbb{R}^N} \rho( \mathbf{x})u^{-\theta}v^{-r}\phi^2dxdy\\
&+\frac{C}{R^2}\int_{ B_{2R}\times B_{(2R)^{1+s}}}\left(u^{-q+1} +
v^{-r+1}\right)dxdy.
  \end{split}
 \end{align}
Fix
\begin{equation}\label{3.6}
 q=\frac{(\theta-1)r}{p-1}-\frac{\theta-p}{p-1}, \quad \mbox{ or equivalently } \quad
1-q=\frac{(\theta-1)(1-r)}{p-1}.
\end{equation}
Let  $r> \frac{p+1}{2}.$  Applying Young's inequality and  using \eqref{3.6},  the  first term on the right hand side of \eqref{3.6}, can be estimated as
\begin{align*}
& \frac{1}{a_1}\int_{\mathbb{R}^N} \rho( \mathbf{x})u^{-q}v^{-p}\phi^2dxdy\\
= & \; \frac{1}{a_1}\int_{\mathbb{R}^N} \rho( \mathbf{x})u^{-\frac{\theta+1}{2}}v^{-\frac{p+1}{2}}u^{-\frac{(\theta-1)r}{p-1}+\frac{\theta-1}{p-1}\left(\frac{1+p}{2}\right)}v^{\frac{1-p}{2}}\phi^2 dxdy \\
= & \;\frac{1}{a_1}\int_{\mathbb{R}^N} \rho( \mathbf{x})u^{-\frac{\theta+1}{2}}v^{-\frac{p+1}{2}}u^{(1-q)\frac{2r-1-p}{2(r-1)}}v^{\frac{1-p}{2}}\phi^2dxdy\\
\leq & \; \frac{2r-1-p}{2(r-1)}\int_{\mathbb{R}^N} \rho( \mathbf{x})u^{-\frac{\theta+1}{2}}v^{-\frac{p+1}{2}}u^{1-q}\phi^2dxdy
 +\frac{p-1}{2(r-1)}a_1^{-\frac{2(r-1)}{p-1}}\int_{\mathbb{R}^N} \rho( \mathbf{x})u^{-\frac{\theta+1}{2}}v^{-\frac{p+1}{2}}v^{1-r}\phi^2dxdy
\\
= & \; \frac{2r-1-p}{2(r-1)}I_1+\frac{p-1}{2(r-1)} a_1^{-\frac{2(r-1)}{p-1}} I_2,
\end{align*}
and similarly
\begin{align*}
  {a_2}^\frac{2r-1-p}{p-1}\int_{\mathbb{R}^N} \rho( \mathbf{x})u^{-\theta}v^{-r}\phi^2dxdy  \leq \frac{p-1}{2(r-1)} I_1
  + \frac{2r-1-p}{2(r-1)}{a_2}^\frac{2(r-1)}{p-1}I_2.
   \end{align*}

Inserting the two above estimates in \eqref{3.5}, we arrive at

\begin{align*}
{a_2}^{\frac{2(r-1)}{p-1}}I_2\leq
\left[\frac{2r-1-p}{2(r-1)}{a_2}^{\frac{2(r-1)}{p-1}}+\frac{p-1}{2(r-1)}{a_1}^{\frac{-2(r-1)}{p-1}}\right]I_2 +
\frac{C}{R^2}\int_{ B_{2R}\times B_{(2R)^{1+s}}}\left(u^{-q+1}
+ v^{-r+1}\right)dxdy.
\end{align*}
As $p\geq\theta>1,$ we have  $r> \frac{p+1}{2}>1,$  and $1-q<0.$ Combining \eqref{3.6} and \eqref{estS}, one obtains
$$u^{1-q} \leq Cv^{1-r}\quad \mbox{and} \quad
u^{-\frac{\theta+1}{2}}v^{-\frac{p+1}{2}}v^{1-r}\geq u^{-\theta} v^{-r}.$$
We get then
\begin{align}
\label{LF}
\frac{p-1}{2(r-1)}\left[(a_1a_2)^{\frac{2(r-1)}{p-1}}-1\right]
\int_{\mathbb{R}^N} \rho( \mathbf{x})u^{-\theta} v^{-r}\phi^2dxdy\leq CR^{-2}a_1^{\frac{2(r-1)}{p-1}}\int_{ B_{2R}\times B_{(2R)^{1+s}}}
v^{1-r}dxdy.
\end{align}
A simple computation find
$$a_1a_2 > 1,$$
is equivalent to
$$16p\theta \left( \frac{\theta-1}{p-1}z^2+\frac{p+\theta-2}{(p-1)^2}z+1\right)>\frac{(\theta-1)^2}{(p-1)^2}z^{4}$$
where we set $z=r-1$, and $z > \frac{p-1}{2}$. That is
\begin{align*}
 L(z):=z^4-\frac{16p\theta(p-1)}{\theta-1}z^2-\frac{16p\theta(p-1)(p+\theta-2)}{(\theta-1)^2}z-\frac{16p\theta(p-1)^2}{(\theta-1)^2}<0.
\end{align*}
Consequently, from \eqref{LF} ,we conclude that
\begin{align*}
 \int_{ B_{R}\times B_{(R)^{1+s}}} \rho( \mathbf{x}) u^{-\theta} v^{-z-1}dxdy  \leq \frac{C}{R^2}\int_{ B_{2R}\times B_{(2R)^{1+s}}} v^{-z}dxdy.
\end{align*}
Furthermore, we can check that $a_{1}a_2>1$ is equivalent to $L(z)<0$, the proof is completed. \qed

\medskip\noindent
\subsection{\bf End of the proof of Theorem \ref{main3}.}
 \medskip

 In this subsection, we use $L^{2}$-estimates for Grushin operator,  and we apply the bootstrap iteration  as in \cite{cow, DP,Hfh}. For the completeness, we present the details.

\smallskip
 Let $ \eta \in C_c^\infty\left(\mathbb{R}^N=\mathbb{R}^{N_1}\times\mathbb{R}^{N_2}, [0,1]\right)$ be  a cut-off function  such that

 \begin{align}
\label{test}
\eta=1 \quad \mbox{on} \quad B_{1}\times B_{1},\quad \mbox{and} \quad  \eta=0 \quad \mbox{outside } \quad B_{2}\times B_{2^{1+s}}.
\end{align}

\smallskip

Let $w$ be a smooth function and let  $\lambda_{s}=\frac{N_{s}}{N_{s}-2}.$ Using Sobolev inequality  \cite{xxy} and integration by parts, we obtain
\begin{align*}
\left(\int_{ B_{1}\times B_{1}} w^{2\lambda_{s}} dxdy\right)^{\frac{1}{2\lambda_{s}}}&\leq \left(\int_{ B_{2}\times B_{2^{1+s}}} (w\eta)^{2\lambda_{s}} dxdy\right)^{\frac{1}{2\lambda_{s}}}\\
&\leq C\left(\int_{ B_{2}\times B_{2^{1+s}}} |\nabla_{s} (w\eta)|^{2} dxdy\right)^{\frac{1}{2}}\\
&\leq C\left[\int_{ B_{2}\times B_{2^{1+s}}} \left(|\nabla_{s} w|^{2}\eta^{2}+w^{2}|\nabla_{s} \eta|^{2}-\frac{w^{2}}{2}\D_{s}(\eta)  \right)dxdy\right]^{\frac{1}{2}}.
\end{align*}
So, we get
\begin{align*}
\left(\int_{ B_{1}\times B_{1}} w^{2\lambda_{s}} dxdy\right)^{\frac{1}{\lambda_{s}}}
\leq C\int_{ B_{2}\times B_{2^{1+s}}} \left(|\nabla_{s} w|^{2}+w^{2}  \right)dxdy.
\end{align*}
By scaling argument, we obtain  readily the estimate

 \begin{align}
\label{test2}
\begin{split}
&\left(\int_{ B_{R}\times B_{R^{1+s}}} w^{2\lambda_{s}} dxdy\right)^{\frac{1}{\lambda_{s}}}\\
&\leq
CR^{N_{s}\big(\frac{1}{\lambda_{s}}-1\big)+2}\int_{ B_{2R}\times B_{(2R)^{1+s}}}|\nabla_{s} w|^{2}dxdy+CR^{N_{s}\big(\frac{1}{\lambda_{s}}-1\big)}\int_{ B_{2R}\times B_{(2R)^{1+s}}}w^{2}dxdy.
\end{split}
\end{align}
\medskip

 Let $(u,v)$ be a  stable solution of \eqref{1.1}, with $1< p \leq \theta$. For $2t_0^-<z_0,$ in what follows, we choose
$$w=v^{-\frac{z_0}{2}}.$$

Let us put $\eta_{R}(x,y)=\eta(\frac{x}{R},\frac{y}{R^{1+s}}),$ where $\eta$ is given in \eqref{test}. By a simple calculation, we obtain readily
 \begin{align}
\label{stbk}
\int_{ B_{R}\times B_{R^{1+s}}} |\nabla_{s} w|^{2} dxdy \leq C \int_{ B_{2R}\times B_{(2R)^{1+s}}}v^{-z_{o}-2}|\nabla_{s} v|^{2}\eta_{R}^{2} dxdy.
\end{align}
 Multiplying $\D_{s} v= \Big(1+\|\mathbf{x}\|^{2(s+1)}\Big)^{\frac{\alpha}{2(s+1)}}u^{-\theta}$ by $v^{-z_{o}-1} \eta_{R}^2$ and  integrating by parts, we derive
\begin{align}\label{3.9}
\begin{split}
& (z_{o}+1) \int_{ B_{2R}\times B_{(2R)^{1+s}}}v^{-z_{o}-2}|\nabla_{s} v|^{2}\eta_{R}^2dxdy\\
&\;= \int_{ B_{2R}\times B_{(2R)^{1+s}}}\Big(1+\|\mathbf{x}\|^{2(s+1)}\Big)^{\frac{\alpha}{2(s+1)}}v^{-z_{o}-1} u^{-\theta}\eta_{R}^2 dxdy\\
&\;+2 \int_{ B_{2R}\times B_{(2R)^{1+s}}}\eta_{R} v^{-z_{o}-1}\nabla_{s}v \cdot \nabla_{s} \eta_{R}  dxdy.
\end{split}
\end{align}
By Young's inequality, we have
\begin{align*}
 &\;  \int_{ B_{2R}\times B_{(2R)^{1+s}}} v^{-z_{0}-1}|\nabla_{s} v||\nabla_{s} \eta_{R}|\eta_{R} dxdy \\
&\; \leq \frac{1}{2} \int_{ B_{2R}\times B_{(2R)^{1+s}}}v^{-z_{0}-2} |\nabla_{s} v|^{2}\eta_{R}^2dxdy +
  \frac{1}{2} \int_{ B_{2R}\times B_{(2R)^{1+s}}}v^{-z_{0}}|\nabla_{s} \eta_{R}| ^2dxdy.
\end{align*}

Substituting this in \eqref{3.9},

\begin{align}\label{9}
\begin{split}
  \int_{ B_{2R}\times B_{(2R)^{1+s}}}v^{-z_{o}-2}|\nabla_{s} v|^{2}\eta_{R}^2dxdy
&\;\leq C \int_{ B_{2R}\times B_{(2R)^{1+s}}}\Big(1+\|\mathbf{x}\|^{2(s+1)}\Big)^{\frac{\alpha}{2(s+1)}}v^{-z_{o}-1} u^{-\theta}\eta_{R}^2 dxdy\\
&\;+C \int_{ B_{2R}\times B_{(2R)^{1+s}}}v^{-z_{0}}|\nabla_{s} \eta_{R}| ^2dxdy.
\end{split}
\end{align}

Combining \eqref{3.1} , \eqref{stbk} and \eqref{9}, one obtains
\begin{align*}
&\;\int_{ B_{R}\times B_{R^{1+s}}} |\nabla_{s} w|^{2} dxdy\\
&\;\leq C\int_{ B_{2R}\times B_{(2R)^{1+s}}}v^{z_{o}-2}|\nabla_{s} v|^{2}\eta_{R}^2dxdy \\
&\;\leq C\int_{ B_{2R}\times B_{(2R)^{1+s}}}\Big(1+\|\mathbf{x}\|^{2(s+1)}\Big)^{\frac{\alpha}{2(s+1)}}v^{-z_{o}-1} u^{-\theta}\eta_{R}^2 dxdy+\frac{C}{R^{2}}\int_{ B_{2R}\times B_{(2R)^{1+s}}}v^{-z_{0}}dxdy
 \\
&\;\leq \frac{C}{R^{2}}\int_{ B_{2R}\times B_{(2R)^{1+s}}}v^{-z_{0}}dxdy.
\end{align*}
In view of estimate \eqref{test2}, we get
\begin{equation}
\label{lnew2}
\left(\int_{ B_{R}\times B_{R^{1+s}}} v^{-z_0\lambda_{s}} dxdy\right)^{\frac{1}{\lambda_{s}}}\leq
CR^{N_{s}\big(\frac{1}{\lambda_{s}}-1\big)}\int_{ B_{2R}\times B_{(2R)^{1+s}}}v^{-z_0}dxdy.
\end{equation}
\medskip

Let $z_0$ be the largest root of the polynomial $L$ given by \eqref{3.1}. We fix a real positive number
 $$q=\frac{p+\theta+2}{2}\quad  \mbox{satisfying} \quad 2t_0^-<q<z_{0},$$
 and let  $m$ be the nonnegative integer such that $$q\lambda_{s}^{m-1}<z_{0}<q\lambda_{s}^{m}.$$
 We construct an increasing geometric sequence
 $$2t_0^-<z_{1}<z_{2}<,....,<z_{m}<z_{0}.$$
 as follows
 $$z_{1}=qk,  \quad z_{2}=qk\lambda_{s},....,z_{m}=qk\lambda_{s}^{m-1},$$
  where $k\in[1,\lambda_{s}],$ will be chosen so that $z_m$ is arbitrarily close to $z_0.$

 Set $ R_n = 2^n R.$ By \eqref{lnew2} and an induction argument, we deduce then
 \begin{align}
\label{lnnnew2}
\begin{split}
\left(\int_{ B_{R}\times B_{R^{1+s}}} v^{-z_m\lambda_{s}} dxdy\right)^{\frac{1}{z_m\lambda_{s}}}&\leq
CR^{N_{s}\big(\frac{1}{z_m\lambda_{s}}-\frac{1}{z_m}\big)}\left(\int_{ B_{1}\times B_{(1)^{1+s}}}v^{-z_m}dxdy\right)^{\frac{1}{z_m}}\\
&\;=CR^{N_{s}\big(\frac{1}{z_m\lambda_{s}}-\frac{1}{z_m}\big)}\left(\int_{ B_{1}\times B_{(1)^{1+s}}}v^{-z_{m-1}\lambda_{s}}dxdy\right)^{\frac{1}{z_{m-1}\lambda_{s}}}\\
&\;\leq CR^{N_{s}\big(\frac{1}{z_m\lambda_{s}}-\frac{1}{z_1}\big)}\left(\int_{ B_{R_{m}}\times B_{(R_{m})^{1+s}}}v^{-z_{1}}dxdy\right)^{\frac{1}{z_{1}}}\\
&\;\leq CR^{N_{s}\big(\frac{1}{z_m\lambda_{s}}-\frac{1}{qk}\big)}\left(\int_{ B_{R_{m}}\times B_{(R_{m})^{1+s}}}v^{-qk}dxdy\right)^{\frac{1}{qk}}
\end{split}
\end{align}
 Furthermore, by H\"older's inequality, there holds
  \begin{align}
\label{lnnnew23}
\begin{split}
&\;\left(\int_{ B_{R_{m}}\times B_{(R_{m})^{1+s}}}v^{-qk}dxdy\right)^{\frac{1}{qk}}\\
&\;\leq
\left[\left(\int_{ B_{R_{m}}\times B_{(R_{m})^{1+s}}}v^{-q\lambda_{s}}dxdy\right)^{\frac{k}{\lambda_{s}}}\left(\int_{ B_{R_{m}}\times B_{(R_{m})^{1+s}}}dxdy\right)^{1-\frac{k}{\lambda_{s}}}\right]^{\frac{1}{qk}}\\
&\;\leq C\left[\left(\int_{ B_{R_{m}}\times B_{(R_{m})^{1+s}}}v^{-q\lambda_{s}}dxdy\right)^{\frac{k}{\lambda_{s}}}CR^{N_{s}\big(1-\frac{k}{\lambda_{s}}\big)}\right]^{\frac{1}{qk}}\\
&\;\leq CR^{N_{s}\big(\frac{1}{kq}-\frac{1}{q\lambda_{s}}\big)}\left(\int_{ B_{R_{m}}\times B_{(R_{m})^{1+s}}}v^{-q\lambda_{s}}dxdy\right)^{\frac{1}{q\lambda_{s}}}\\
&\;\leq CR^{N_{s}\big(\frac{1}{kq}-\frac{1}{q\lambda_{s}}\big)}R^{N_{s}\big(\frac{1}{q\lambda_{s}}-\frac{1}{q}\big)}\left(\int_{ B_{R_{m}}\times B_{(R_{m})^{1+s}}}v^{-q}dxdy\right)^{\frac{1}{q}}.
\end{split}
\end{align}
 Combining the last tow inequalities, we obtain
 \begin{align}
 \label{lnew2j}
 \left(\int_{ B_{R}\times B_{R^{1+s}}} v^{-z_m\lambda_{s}} dxdy\right)^{\frac{1}{z_m\lambda_{s}}}\leq
CR^{N_{s}\big(\frac{1}{z_m\lambda_{s}}-\frac{1}{q}\big)}\left(\int_{ B_{R_{m}}\times B_{(R_{m})^{1+s}}}v^{-q}dxdy\right)^{\frac{1}{q}},
\end{align}

We deduce from \eqref{VV} that

\begin{align}
 \label{lmbv}
 \left(\int_{ B_{R}\times B_{R^{1+s}}} v^{z_m\lambda_{s}} dxdy\right)^{\frac{1}{z_m\lambda_{s}}}\leq
CR^{\frac{N_{s}}{z_m\lambda_{s}}-\frac{2+\alpha}{q}}.
\end{align}
 Recall that $\lambda_{s}=\frac{N_{s}}{N_{s}-2}.$ Suppose now
$$N_{s}<2+2\left(\frac{2+\alpha}{p+\theta +2}\right)z_0,$$
 we can chose $k\in[1,\lambda_{s}],$ such that $z_{m}$ is sufficiently  close to  $z_{0}$ satisfying $$N_{s}-2-2\left(\frac{2+\alpha}{p+\theta +2}\right)z_m<0.$$
  Then, it implies from \eqref{lmbv} that $\|v\|_{L^{z_m\lambda_{s}}(\R^N)} = 0$ as $R\rightarrow\infty,$ i.e., $v \equiv 0$ in $\R^N$ This is a contraction.
  \smallskip

Therefore, we get the desired result, the equation \eqref{1.1} has no
stable solution if $$N_{s}<2\Big[1+(2+\alpha)x_0\Big], \quad \mbox{where} \quad  x_0=\frac{1}{p+\theta+2}z_0.$$

\medskip\noindent
\textbf{Proof of Corollary \ref{main2}.} Let $u$ be a stable solution of equation \eqref{D}, then $v = u$ verify the system \eqref{1.1} with $p = \theta$.  Moreover, we have $$t_0^{\pm}= p\pm\sqrt{p^2+p},$$ and
\begin{align*}
  L(z) = z^4 -16p^2z^2 +
32p^2z-16p^2=(z^2+4p(z+1))(z-2t_0^-)(z-2t_0^+).
\end{align*}

As $t_0^+ > p > 1,$ it follows that $2t_0^+$ is the largest root of $L$ as $t_0^+ > p > 1$. Therefore
$$x_0 =\frac{2t_0^+}{p+1}= \frac{2p+2\sqrt{p^2+p}}{p+1}$$ is the largest root of $H.$ Then, applying Theorem \ref{main3}, the result follows immediately. \qed

\end{document}